# A Distributed Power Routing Method between Regional Markets based on Bellman-Ford Algorithm

Xue Tai, Hongbin Sun, Senior Member, IEEE, Qinglai Guo, Senior *Member, IEEE* and Zihao Li

*Abstract*—With the deregulation of power industry, large power users can buy electricity directly from power producers. For trades between different regional markets, it is necessary to research on the method to find the cheapest route for power exchanging. Moreover, as there exists no coordinator above regional markets, the power routing problem must be solved in a distributed way. Most previous research are based on power electronic transformers. Some research proposed EMS of sub-networks as energy routers, but failed to consider the power flow constraints or did not converge in some specific cases. Our method is based on the EMS of multiple markets, which saves cost by eliminating expensive electronic transformers and has high efficiency in systems with numerous nodes by considering the regional markets as generalized nodes. We take transmission costs and congestion costs into consideration, and uses Bellman-Ford algorithm to determine the cheapest route for power transmission between different regional markets and make convergence analysis about our algorithm. This method proves effective and convergent when there is no loop among the markets. A four-market loop-free case was used to illustrate the efficiency of the method.

## I. INTRODUCTION

WITH the deregulation of power industry, large power users can purchase electricity directly from energy producers.[1] Because of different electricity prices in different regional markets, sometimes energy trades may occur between different markets in order to make full use of cheap renewable power such as hydroelectricity and wind power. Therefore, it is important for market participants to find the cheapest way to transmit the power to reduce costs and avoid congestion.

Conventional power routing method involves power electronic transformer to change the route of power flow [2]-[6]. [2] proposed a dynamic power routing scheme for hybrid ac/dc microgrids operating in islanded mode which utilized the interlinking converters between the ac and dc sides. [3] focused on the issues of access management and consuming in energy routing network within the low-voltage distribution grid and proposed the routing matrix of energy router to control the power flow. [4] proposed the idea of power router which is capable of direct the power towards a selected path regardless of the resistance of the power line. It also addressed on the problem of communication time delay. [5] proposed an on-demand electric power supply architecture in home based on quality-aware routing. [6] used Dijkstra algorithm to search the [1]shortest path to implement power routing and optimize system operation. All the research above relies on power electronic transformers to change the power flow, in which case every node involved in energy transactions needs to be equipped with power electronic devices, which is relatively expensive at the present stage.

The power grid network has a lot in common with the Internet, both contains millions of nodes and links between them. The transmission of power is analogous to the transmission of data in the Internet. Therefore, some routing methods in the Internet can be applied to power transactions in the power system. In the field of Internet, routing algorithms include static routing and dynamic routing.[7] Static routing is only suitable for networks that do not change over time. In the power system, the topology of networks usually changes due to different operation modes or unexpected accidents. Therefore, we should applied dynamic routing algorithm in finding the cheapest route for a transaction. Widely-used dynamic algorithm include distance vector algorithms, link-state algorithms, optimized link state routing algorithm and path vector protocol. Distance vector algorithms use the Bellman-Ford algorithm, which is able to find the shortest route from two certain nodes in a decentralized way.[8] Link-state algorithms requires all nodes to form a graphical map of the network, and calculate the shortest path to every other node, which consumes two much computing power compared with distance vector algorithms, since the latter only needs to find the shortest path between two certain nodes.[9] Optimized link state routing algorithm is similar to link state algorithm but is usually used in wireless networks.[10] Path vector algorithm is often used with Bellman-Ford algorithm to maintain and update the path information.[11] [12][13]converted the traditional OPF problem into minimum cost flow problem and applied shortest

This work is supported in part by The National Key Research and Development Program of China (Basic Research Class 2017YFB0903000) and in part by funding of State Grid Corporation of China (Basic Theories and Methods of Analysis and Control of the Cyber Physical Systems for Power Grid).

X. Tai, H. Sun , Q. Guo and Z. Li are with the State Key Laboratory of Power Systems, Department of Electrical Engineering, Tsinghua University, Beijing 10084, China. (Corresponding author: Hongbin Sun, shb@tsinghua.edu.cn)

path algorithm to manage power in multi-agent system based active distribution networks. However, according to the count-to-infinity problem proposed in [8], their algorithm did not converge in some specific cases. Therefore, based on the work above, we use Bellman-Ford algorithm to solve the route searching problem in power transaction, and made some adjustments to the algorithm to ensure its convergence.

[14] presented a distributed control and coordination architecture whose key idea is to distribute the intelligence into the periphery of the grid. It proposed the concept of clusters which are collection of electric grid system resources. Two main cluster types are a resource cluster and a balanced cluster. A resource cluster is a collection of resources that supply a service, while a balanced cluster is a collection of resources that largely achieve internal power balance. Any mismatch between power consumption and production is met by contracted imports and exports from other clusters. In short, any system with the ability to manage the energy can be regarded as a cluster. In our case, a regional market can be considered as a cluster. It has an intelligent periphery in that it can not only rearrange the output of generators inside the market when a transaction occurs, but can also communicate with other markets to find the cheapest route to transmit the trading energy to another market and make adaptions to relieve congestion in the meantime. Contrary to traditional concept of power electronic transformers as energy router, a regional market with its EMS can serve as an energy router because of its ability to control the power flow to some extent by adjusting the output of generators and discover the lowest-cost route for power transmission.

[15] –[16] used the concept of clusters and proposed a novel cloud-based approach to deal with the optimal power routing problem in clusters of DC microgrids. However, they did not consider the power flow constraints in their optimization problem, and the method they proposed cannot solve the routing problem in a decentralized way. Therefore, new approaches should be studied to solve the power routing problem.

We consider regional markets with EMS as clusters which serve as energy routers. We propose a power routing approach based on Bellman-Ford algorithm in the article. The remainder of this paper is organized as follows: section II introduce the mathematical model of the optimized routing problem. Section III describes the Bellman-Ford algorithm and applies it to the power routing problem. Section IV uses a four-market case to illustrate the routing algorithm proposed. Section V makes a discussion about the algorithm.

## II. MATHEMATIC MODEL

All regional markets that constitute the whole power market are modeled as an undirected graph G (V, E). V stands for vertex, which is the regional market. E stands for edge, which means the transmission lines between the markets.

Transmission costs contain two parts: transmission fee and congestion cost.

### A. Transmission Fee

A transaction that passes through a particular regional market should be charged a cost linear to the power, which can be written as follows:

$$C_{ti} = \text{Price}_i \cdot P \quad (1)$$

$C_{ti}$ is the transmission cost charged by regional market i, when a transaction passes through market i. $\text{Price}_i$ means the money regional market i charged for every megawatt's power that passes through it. P means the amount of power passing through market i.

A transaction that passes through a power line should be charged a cost linear to the power, which can be written as follows:

$$C_{tij} = \text{Price}_{ij} \cdot P \quad (2)$$

$C_{tij}$ stands for transmission cost through the power line between market i and market j. $\text{Price}_{ij}$ means the power line charged for every megawatt's power that passes through it. P means the amount of power passing through the power line.

### B. Congestion Fee

A transaction that causes congestion of a regional market should be charged congestion fee. Congestion fee is the extra cost of all the generators in a market caused by a transaction, which is calculated as below:

$$C_{ci} = \min \sum C_t(G_i) - \min \sum C(G_i) \quad (3)$$

$C_{ci}$ stands for congestion fee charged by regional market i. $\sum C_t(G_i)$ means the total cost of generators in market i when a transaction passes through the market. $\sum C(G_i)$ denotes the total cost of generators when there is no transaction passing through the market.

Markets that on the route of the transaction can be divided into three categories: source market, intermediate market and target market. The source market means the market where the seller of the transaction exists, and the transaction power flows out of it. An intermediate market is the market where neither seller nor buyer of the transaction exists, and the power just passes through it. The source market is the market where the buyer of the transaction exists, and transaction power flows into it. When calculating the total cost of generators when there is a transaction, the three categories of markets have different expressions.

For source markets, the total cost of generators can be calculated as follows:

$$\min \sum C_t(G_i)$$
$$s.t. \quad P_{seller} \geq P_{tr}$$
$$P_{out} = P_{outload} + P_{tr} \quad (4)$$
$$P_i = B_i \theta_i$$
$$|P_{jk}| \leq P_{jk\max} \quad \forall j,k \in i$$

$P_{tr}$ means the power of the transaction. $P_{seller}$ means the injection of the seller node in the transaction. $P_{out}$ means the adapted injection of the node through which the transaction power flows out of the market. $P_{outload}$ is the original load of the node through which the transaction flows out of the market. $B_i$

is the B matrix for a DC power flow in market i. $\theta_i$ is the angle of each bus in market i. $P_{jk}$ is the power flow through the power line between node j and node k in market i. $P_{jkmax}$ is the maximum power flow that the power line can hold.

For intermediate markets, the total cost of generators can be calculated as follows:

$$\min \sum C_t(G_i)$$
$$s.t. \ P_{in} = P_{inload} - P_{tr}$$
$$P_{out} = P_{outload} + P_{tr} \quad (5)$$
$$P_i = B_i \theta_i$$
$$|P_{jk}| \leq P_{jk\max} \quad \forall j,k \in i$$

$P_{in}$ means the adapted injection of the node from which the transaction power flows into the market. $P_{inload}$ means the original load of the node through which the transaction flows into the market.

For target markets, the total cost of generators can be calculated as follows:

$$\min \sum C_t(G_i)$$
$$s.t. \ P_{in} = P_{inload} - P_{tr}$$
$$P_{buyer} = P_{buyerload} + P_{tr} \quad (6)$$
$$P_i = B_i \theta_i$$
$$|P_{jk}| \leq P_{jk\max} \quad \forall j,k \in i$$

$P_{buyer}$ means the adapted injection of the buyer node of the transaction. $P_{buyerload}$ means the original load of the buyer node.
(3) optimization problem

## C. Optimization Problem

The power routing problem is to find a route for power to be delivered from the producer to the consumer. Using DC power flow, the optimization problem can be written in the following form:

$$\min_R \sum_{i \in R}(C_{ti} + C_{ci})$$
$$s.t. \ P = B\theta \quad (7)$$
$$|P_{ij}| \leq P_{ij\max} \quad \forall i,j \in V$$

Where P stands for the power injection of every node in the power network. B is the B matrix for a DC power flow. $\theta$ is the angle of each bus. $P_{ij}$ is the power flow of the power line between bus i and j. $P_{ijmax}$ is the maximum power flow that the power line can hold.

To solve the optimization problem in a decentralized way, we used the Bellman-Ford algorithm to find the cheapest route for power transmission, which will be elaborated on in the next section.

## III. BELLMAN-FORD ALGORITHM

Bellman-Ford algorithm is an algorithm that computes shortest paths from a single source vertex to other vertices in a weighted digraph. In the algorithm, the approximation to the correct distance to a vertex is gradually replaced by more accurate values through iterations until eventually reaching the optimum solution. The Bellman-Ford algorithm relaxes all the edges, and in loop-free networks, by iterating at most $|V|-1$ times, an optimum solution can be found, where $|V|$ is the number of vertices in the graph.[17]

The iteration times depends on the relative iterating order of each node in every step of iteration. In the loop-free system, when in every step, the iterating order happens to be from the power seller to the power consumer, the iteration times is two. However, if the iterating order happens to be from power consumer to the power seller. The iteration time is $|V|-1$, which is the worst case. Therefore, the converging time of the algorithm is linear to the number of markets, which is acceptable in practice.

However, when the topology of the graph changes, there may be count-to-infinity problem as the node does not record the information of the route to a certain node[8]. To avoid this problem, we adapt the Bellman-Ford problem by having the nodes record and update their route information during every iteration.

In the power routing problem described above, the regional markets can be regarded as vertices. The weight of edges between the vertices is the transmission cost between the markets. The weight of edges can be written in the following form:

$$D_{ij} = C_{tj} + C_{tij} + C_{ci}$$

which is the cost of power flowing from market i to market j. The power routing problem is transformed into finding the shortest path from two vertices.

The step of Bellman-Ford algorithm applied to the power routing problem is as follows, assuming the transaction is from market i to market j:
1) Each node calculate its distance to node i. If there is an edge between node k to node i, its distance is $D_{ik}$, otherwise the distance is infinite.
2) Each node communicate with its neighbors. If for node k, there exists a neighbor node m, and $D_{ik} > D_{im} + D_{mk}$, set the distance of node k $D_{ik}$ to $D_{im}+D_{mk}$. Meanwhile, node m informs node j with the information about its path to node i so that node j can update its cheapest route to node i.
3) Check if any node's distance to node i has changed during step 2. If there is some change, go to step 2. Otherwise, continue.
4) Node j gets its shortest distance $D_{ij}$ and its route to node i. It communicate with node i to settle the transaction and make payments to the markets the power passes through.

It is worth noting that, in some cases, due to congestion, there may be no validate route, which will be elaborated on in the following section.

When there are more than one transaction, the route searching problem can be solved in the first-in-first-serve order. The calculation of the extra cost of generators caused by the current transaction is always based on the power flow influenced by the previous transactions.

## IV. CASE STUDY

The case consists of four markets. Each market is a 9-node system. The four markets constitute a loop-free graph, which is shown as below. The transmission cost for each power line between the four markets are all set to 1$/MWh. The transmission cost for each market is set to 1$/MWh.

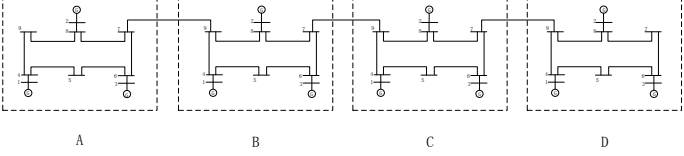

Fig. 1. A loop-free system made up of four regional markets

### A. No Congestion

There is no transaction between the four markets in the beginning. When there is a transaction from market A to market D, Bellman-Ford method is adopted to find the shortest route from A to D. The transaction amount is 100MW/h. The iteration process is shown in the following table, and the shortest route is A-B-C-D.

TABLE I
ITERATION PROCESS OF A SINGLE TRANSACTION ADDING ON THE SYSTEM

| Iteration # | $D_{AB}$ | $D_{AC}$ | $D_{AD}$ |
|---|---|---|---|
| 1 | 2949.02 | ∞ | ∞ |
| 2 | 2949.02 | 3149.02 | 3349.02 |

The total cost of generators of every market is shown as below.

TABLE III
COST OF GENERATORS IN EACH MARKET WITHOUT CONGESTION

| Market | Transaction cost |
|---|---|
| A | 7965.05 $/h |
| B | 5216.03 $/h |
| C | 5216.03 $/h |
| D | 5216.03 $/h |

### B. Occurrence of Congestion

When we changed the limit power of line 3-6 in the four markets to 100MW, because of the congestion the transaction causes, the cost of generators in market A, B and C will increase, which is shown as below.

TABLE III
COST OF GENERATORS IN EACH MARKET WITH CONGESTION

| Market | Transaction cost |
|---|---|
| A | 8103.13 $/h |
| B | 5505.70 $/h |
| C | 5505.70 $/h |
| D | 5216.03 $/h |

When we further change the limit power of line in the four markets, because of congestion, the power cannot flow through market B and market C. Therefore, the transaction cannot be settled successfully.

### C. Multiple Transactions

When there are two transactions adding on the system, Bellman-Ford algorithm is conducted successively. The second transaction is from market C to market B. The transaction amount is 150MW/h. The iteration process for the second transaction is shown as below.

TABLE IV
ITERATION PROCESS OF THE SECOND TRANSACTION ADDING ON THE SYSTEM

| Iteration # | $D_{CA}$ | $D_{CB}$ | $D_{CD}$ |
|---|---|---|---|
| 1 | ∞ | 4625.3 | 4625.3 |
| 2 | 4825.3 | 4625.3 | 4625.3 |

The total cost of generators of every market is shown as below.

TABLE V
COST OF GENRATORS WITH MULTIPLE TRANSACTIONS

| Market | Transaction cost |
|---|---|
| A | 7965.05 $/h |
| B | 5216.03 $/h |
| C | 5216.03 $/h |
| D | 9641.33 $/h |

## V. DISCUSSION

From the case study, it can be shown that when there is no congestion, the transmission cost is only made up of wheeling cost charged by the regional market that the transaction passed through. When congestion occurs, transmission cost include both wheeling cost and congestion cost, which leads to the increase of transaction cost.

The iteration times is small for each transaction. Therefore, the shortest route can be searched efficiently. When the limit of power line is too small, there may be no valid solution to the route searching problem. The transaction will be denied.

When there are more than one transaction, the problem can be solved in the order that the transaction is formed. It should be noticed that the calculation of transaction cost and congestion cost needed to be based on the previous transactions that have been settled. It can be seen that in tab.3, the transaction cost of market A is the same as that in tab.4, which means that the transmission cost of the second transaction is based on the first transaction.

Meanwhile, using the Bellman-Ford algorithm, the route-searching problem can be solved in a distributed way, and no control center is needed to coordinate the transactions, which provides a cheap and transparent method for transactions between regional markets.

## VI. CONCLUSION

In this paper, we proposed a novel power routing algorithm based on Bellman-Ford algorithm to solve the problem in which two identities from different regional markets need to find a route for their transaction. When there is no loop among the regional markets, the algorithm can converge in the time linear to the number of markets, which is acceptable in practice. Case study shows that the algorithm is efficient and of high performance in finding a route of a transaction. Further research should be carried out considering systems with loops.